\documentclass[12pt, twoside]{article}
\usepackage{amsmath,amsthm}
\usepackage{amsfonts}
\usepackage{amssymb,latexsym}
\usepackage{enumerate}
\usepackage{cases}
\newtheorem{theorem}{Theorem}[section]

\newtheorem{proposition}[theorem]{Proposition}
\newtheorem{lemma}[theorem]{Lemma}
\newtheorem{definition}[theorem]{Definition}

\newtheorem{corollary}[theorem]{Corollary}

\def\ind{1{\hskip -3 pt}\hbox{\textsc{I}}}
\def\n{\noindent}

\def\Om{\Omega}
\def\E{\mathcal E}

\def\ve{\varepsilon}
\def\va{\varphi}
\def\C{\mathbb{C}^n}

\def\O{\widetilde{\Omega}}
 
\def\u{\tilde{u}} 
\def\w{\hat{w}}
 
\def\F{\mathcal F}
\DeclareMathOperator{\loc}{loc}
\numberwithin{equation}{section}
\usepackage{scrextend}
\changefontsizes[16pt]{13pt}
\begin{document}
\setlength{\baselineskip}{18truept}
\pagestyle{myheadings}
\markboth{ N. V. Phu }{ maximal subextension and stability on $m$-capacity }
\title {Maximal subextension and stability on $m$-capacity of maximal subextension of $m$-subharmonic functions with given boundary values.  }
\author{
 Nguyen Van Phu
\\ Faculty of Natural Sciences, Electric Power University,\\ Hanoi,Vietnam.
\\E-mail: phunv@epu.edu.vn}

\date{}
\maketitle

\renewcommand{\thefootnote}{}



\renewcommand{\thefootnote}{\arabic{footnote}}
\setcounter{footnote}{0}

\begin{abstract}
\n  In this paper, we study maximal subextension of $m$-subharmonic functions with given boundary values. We also prove stability on $m$-capacity of maximal subextension of $m$-subharmonic functions with given boundary values.
\end{abstract}

\section{Introduction}
Let $\Omega\subset \tilde{\Omega}$ be bounded domains in $\mathbb{C}^{n}$ and $u\in PSH(\Omega).$ A function $\tilde{u}\in PSH(\tilde{\Omega})$ is said to be a subextension of $u$ from $\Om$ to $\O$ if $\tilde{u}\leq u$ on $\Omega.$
Subextension in the class $\mathcal{F}(\Omega),$ where $\Omega$ is a bounded hyperconvex domain in $\mathbb{C}^{n}$ has been studied by Cegrell and Zeriahi in \cite{CZ03}. In details, the authors proved that if $\Omega\Subset \tilde{\Omega}$ are bounded hyperconvex domains in $\mathbb{C}^{n}$ and $u\in\mathcal{F}(\Omega)$ then there exists $\tilde{u}\in\mathcal{F}(\tilde{\Omega})$ such that $\tilde{u}\leq u$ on $\Omega$ and $\int_{\tilde{\Omega}}(dd^{c}\tilde{u})^{n}\leq \int_{\Omega}(dd^{c}u)^{n}.$ For the class $\E_p(\Om),p>0,$ the subextension problem was investigated by P. H. Hiep \cite{H08}, who proved that if $\Om\subset\O$ are hyperconvex domains and $u\in\E_p(\Om),p>0,$ then there exists a function $\u\in\E_p(\O)$ such that $\u\leq u$ on $\Om$ and $\int_{\O}(-\u)^p(dd^c\u)^n\leq \int_{\Om}(-u)^p(dd^cu)^n.$ The subextension problem concerning boundary values in bounded hyperconvex domains was studied by  Czy\.{z} and Hed \cite{CzH08} and after that by Hai and Hong \cite{HH14}. In 2015, Chinh \cite{Ch15} introduced and investigated the Cegrell classes of $m$-subharmonic functions which are more generalized than the Cegrell classes of plurisubharmonic functions. Subextension of $m$-subharmonic functions in the class $\F_m(\Om)$ in bounded $m$-hyperconvex domains was studied by Hai and Dung in \cite{HD20}. 
 After that, Phu and Dieu \cite{PD23} studied subextension in the class $\E_{m,\chi}(\Om).$ Recently, N. V. Phu \cite{Pjmaa} studied subextension of $m$-subharmonic functions with given boundary values $\mathcal{F}_m(f,\Omega).$ In detail, Theorem 4.4 in \cite{Pjmaa} proved that: Let $\Om\subset \O$ are bounded $m-$hyperconvex domains in $\C.$ Assume that $f\in\E_m(\Om)$ and $g\in\E_m(\O)\cap MSH_m(\O)$ satisfying $f\geq g$ on $\Om.$  Then for every $u\in\F_m(f,\Om)$ satisfying $\int_{\Om}H_m(u)<\infty,$  there exists a subextension $\phi$ of $u$ to $\O$ satifying $\phi\in \mathcal{F}_{m}(g,\widetilde{\Omega})$ and $H_m(\phi) \leq\ind_{\Om} H_m(u) $ on $\O.$\\
 \n On the other hand, finding the conditions for a sequence of $m$-subharmonic functions $[u_j]$ such that the sequence $[H_m(u_j)]$ converges in the weak* topology is one of the most important problems. It is known that the weak convergence of a sequence of $m$-subharmonic functions $[u_j]$ does not give the convergence of the associated Hessian measure $H_m(u_j)$ in general. An essential tool in the study of continuity of Hessian operator is the convergence in $m$-capacity.  A result of Chinh $\cite{Ch15}$ asserts that the convergence in $m$-capacity for bounded sequence of $m$-subharmonic functions is sufficient to obtain the convergence of $[H_m(u_j)].$ Later, Hung and Phu $\cite{HP17}$ proved a similar results in the case $[u_j]\subset\E_m(\Om)$ is bounded from below by a function in the Cegrell's class $\E_m(\Om)$.   
 Recently, Hbil and Zaway $\cite{HZ23}$ gave equivalent statements of convergence in $m$-capacity for sequence of $m$-subharmonic functions with give boundary values when it converges almost everywhere to an $m$-subharmonic function.\\
\n Motivated by the above results, in the literature, firstly, we would like to carry out the subextension with given boundary values for the class of $m$-subharmonic functions. We will improve the result in \cite{Pjmaa} with the general condition of boundary value function $g$. In detail, while $g\in\E_m(\O)\cap MSH_m(\O)$ in Theorem 4.4 in \cite{Pjmaa}, in this paper, we only consider condition $g\in\E_m(\O).$ Secondly, we study stability on $m$-capacity of maximal subextension with given boundary values of $m$-subharmonic functions. Indetail, we will prove that the convergence in $m$-capacity of a sequence of $m$-subharmonic functions with given boundary value implies the convergence in $m$-capacity of their corresponding maximal subextensions with given boundary values.\\
\n The paper is organized as follows. Besides the Introduction, the paper has other four sections. In Section 2 we recall the definitions and results concerning the $m$-subharmonic functions which were introduced and investigated intensively in the recent years by many authors (see \cite{Bl1}, \cite{SA12}, \cite{Ch12}, \cite{T19} and \cite{Ch15}). In Section 3, we present subextension of $m$-subharmonic functions with boundary values. In Section 4, we prove the convergence in $m$-capacity of $m$-subharmonic functions with given boundary values implies the convergence in $m$-capacity of the corresponding maximal subextensions with given boundary values.

\section{Preliminaries}
 The elements of the theory of $m$-subharmonic functions and the complex $m$-Hessian operator can be found e.g. in \cite{Bl1}, \cite{Ch12}, \cite{Ch15}, \cite{SA12} and \cite{T19}. A summary of the properties required for this paper can be found in Preliminaries Section (from subsection 2.1 to subsection 2.6) in \cite{Pjmaa}.\\
 
\begin{definition}
{\rm
 Let $\Om\subset\O$ are $m$-hyperconvex domains and $u\in SH_m(\Om).$ Assume that $f\in\E_m(\Om),g\in\E_m(\O)$ are such that $f\geq g$ on $\Om.$  Let $u\in\F_m(f,\Om),$ consider the function $M_{u,g}$  defined by 
$$M_{u,g}=\big(\sup\{v\in SH^{-}_{m}(\tilde{\Omega}): v\leq u \,\,\text{on}\,\, \Omega \,\,\text{and}\,\, v\leq g \,\,\text{on}\,\, \widetilde{\Omega} \}\big)^*.$$ 
If $M_{u,g}$ is existent in the above definition then $M_{u,g}$ is said to be the maximal subextension of $u$ from $\Om$ to $\O$ with boundary value $g.$}
\end{definition}
\n  When $g\equiv 0$, to make it convenient for presentation, we  use symbol $M_u$ instead of $M_{u,0}.$
\begin{definition}
	{\rm  Assume that $\mathcal{K}\in\{\mathcal{E}_m^0,\mathcal{E}_m,\mathcal{F}_m,\mathcal{N}_m \}$ and $f\in\mathcal{E}_m(\Omega),$ we say that a $m-$subharmonic $u$ defined on $\Omega$ belongs to $\mathcal{K}(f)=\mathcal{K}(f,\Omega)$ if there exists a function $\varphi\in\mathcal{K}$ such that $f\geq u\geq \varphi +f.$
	}
\end{definition}
\n Moreover, if $\varphi\in\mathcal{N}_m^a(\Omega)$ (resp., $\varphi\in\mathcal{F}_m^a(\Omega)$ ) then we say that $u\in\mathcal{N}_m^a(f,\Omega)$ (resp., $u\in\mathcal{F}_m^a(f,\Omega).$)\\
\n We recall some results that will frequently be used in this paper.
\begin{theorem}[Theorem 3.8 in $\cite{HP17}$]\label{th3.8HP17}
	Let	$u_j,v_j,w\in\E_m(\Om)$ be such that $u_j,v_j\geq w,\forall j\geq 1.$ Assume that $|u_j-v_j|\to 0$ in $Cap_m-$ capacity. Then we have $hH_m(u_j)-H_m(v_j)\to 0$ weakly as $j\to+\infty$ for all $h\in SH_m\cap L^{\infty}_{\loc}(\Om).$
\end{theorem}
\n Lemma 2.9 in $\cite{T19}$ stated that if $[u_j]$ is a monotone sequence of $m$-subharmonic functions converging to an $m$-subharmonic function $u$ then we have $u_j\to u$ in $Cap_m$ as $j\to\infty.$ Thus, we have the following Corollary.
\begin{corollary}\label{co3.8HP17}
If $[u_j]\subset\E_m(\Om)$ is a monotone sequence of functions converging to a function $u \in\E_m(\Om)$ then the sequence of measures $H_m(u_j)\to H_m(u)$ weakly as $j\to\infty.$
\end{corollary}
\begin{proposition}[Proposition 2.9 in \cite{Pjmaa}]\label{pro5.2HP17}
\n  Assume that $u,v,u_k\in\E_m(\Om), k=1,\cdots,m-1$ with $u\geq v $ on $\Om$ and  $T=dd^cu_1\wedge\cdots\wedge dd^cu_{m-1}\wedge\beta^{n-m}.$ Then we have
$$\ind_{\{u=-\infty\}}dd^cu\wedge T\leq \ind_{\{v=-\infty\}}dd^cv\wedge T.$$
In particular, if $u,v\in\E_m(\Om)$ are such that $u\geq v$ then for every $m-$polar set $A\subset\Om$ we have
$$ \int_{A}H_m(u)\leq\int_{A}H_m(v).$$
\end{proposition}

\begin{lemma}[Lemma 2.7 in $\cite{Gasmi}$]\label{lm2.7Gasmi}
Let $\varphi\in SH_m^-(\Om)$ and $u,v\in\mathcal{N}_m(f)$ are such that $u\leq v$ and $\int_{\Omega}-\varphi H_m(u)<+\infty.$ Then the following inequality holds
$$\int_{\Om}-\varphi H_m(v)\leq\int_{\Om}-\varphi H_m(u).$$
\end{lemma}


\n We have the following usefull Lemma. 
\begin{lemma}\label{lm4.6PDjmaa}
	Let $\mu$ be a positive measure on $\Om$ which vanishes on all $m-$polar sets and $\mu(\Om)<\infty.$ Let ${u_j}\in SH_m^-(\Om)$ be a sequence satifying the following conditions:\\
	\n (i) $\sup\limits_{j\geq 1}\int_{\Om}-u_jd\mu<+\infty$\\
	\n (ii) $u_j\to u\in SH^-_m(\Om)$ a.e. $dV_{2n}.$\\
	Then we have $\lim\limits_{j\to\infty}\int_{\Om}u_jd\mu=\int_{\Om}ud\mu.$
\end{lemma}
\begin{proof}
	See Step 1 in the proof of Lemma 4.6 in $\cite{PD23jmaa}.$
\end{proof}

\section{Maximal subextension with given boundary values of $m$-subharmonic functions}
\n In this Section, we will prove the following Theorem. This is an extensive version of Theorem 4.4 in \cite{Pjmaa}.
\begin{theorem}\label{t3.3}
Let $\Om\subset\O$ are bounded $m-$ hyperconvex domains in $\mathbb{C}^n.$ Assume that $f\in\E_m(\Om),g\in\E_m(\O)$ are such that $f\geq g$ on $\Om.$\\
{\bf (Part A)} Let $u\in\F_m(f,\Om)$ satisfying $\int_{\Om}H_m(u)<+\infty.$ Then the maximal subextension of $u$ from $\Om$ to $\O$ with boundary value $g$ exists and the following statements hold true:\\
\n (i) $M_{u,g}\in \mathcal{F}_{m}(g,\widetilde{\Omega}).$\\
\n (ii)   $H_m(M_{u,g})=0$ on $\{M_{u,g}<\min(u,g)\}\cap\Om.$\\ 
\n (iii) $H_m(M_{u,g})\leq \ind_{\Om}H_m(u) +H_m(g)$ on $\O.$\\
{\bf (Part B)} Let $u\in\F^a_m(f,\Om)$ satisfying $u\geq g$ on $\Om\smallsetminus K$ for some compact subset $K$ of $\Om.$  Then $M_{u,g}$ exists and the following statements hold true:\\
\n (i) $M_{u,g}\in \mathcal{F}^a_{m}(g,\widetilde{\Omega}).$\\
\n (ii)   $H_m(M_{u,g})=0$ on $\big(\{M_{u,g}<\min(u,g)\}\cap\Om\big)\cup \big((\O\smallsetminus K)\cap\{-\infty<M_{u,g}<g\}\big).$\\ 
\n (iii) $H_m(M_{u,g})\leq \ind_{K\cap\{M_{u,g}=u\}}H_m(u) +H_m(g)$ on $\O.$
\end{theorem}

\begin{proof}
\n {\bf The proof of Part (A)}\\	
 \n (i) By the definition of $u\in\F_m(f,\Om),$ there exists a function $\varphi\in\F_m(\Om)$ such that
 $$f\geq u\geq f+\varphi.$$	
 We put
 $$\hat{\varphi}=\big(\sup\{v\in SH^{-}_{m}(\tilde{\Omega}):  v\leq \varphi\,\,\text{on}\,\, \Om \}\big)^*.$$
 By repeating the argument as in Theorem 3.3 of $\cite{PD23}$ (when $\chi\equiv -1$), we have $\hat{\varphi}\in \F_m(\O),\hat{\varphi}\leq\varphi$ on $\Om$ and $H_m(\hat{\varphi})\leq \ind_{\Om}H_m(\va).$\\
 Now, we have $ g+\hat{\varphi}\leq f+\varphi\leq u$ on $\Om$. Obviously, $g+\hat{\varphi}\leq g$ on $\O.$ Thus, $g+\hat{\varphi}$ belongs to the class of functions which is in the definition of $M_{u,g}.$ So we have $g+\hat{\varphi}\leq M_{u,g}$ on $\O.$ That means we have $M_{u,g}\in \mathcal{F}_{m}(g,\widetilde{\Omega}).$ This implies that we have (i).\\
 \n (ii)  We see that
$$	\{M_{u,g}<\min(u,g)\}\cap\Om = \bigcup_{s\in\mathbb{Q}^-}\big(\{M_{u,g}<s<\min(u,g)\}\cap\Om\big)$$
Thus, it suffices to prove that
$$H_m(M_{u,g})=0\,\,\text{on}
\,\, \{M_{u,g}<s<\min(u,g)\}\cap\Om.$$
 For this, by Theorem 3.1 in $\cite{Ch15}$, there exists a sequence $u_j\in\mathcal{E}_m^0(\Om)\cap C(\overline{\Om})$ such that $u_j$ decreases to $u$ on $\Om,$ supp$H_m(u_j)\Subset\Om$ and a sequence $g_j\in\mathcal{E}_m^0(\O)\cap C(\overline{\O})$ such that $g_j$ decreases to $g$ on $\O,$ supp$H_m(g_j)\Subset\O.$
 Consider function $\w_j$ defined by
$$\w_j=\big(\sup\{v\in SH^{-}_{m}(\tilde{\Omega}): v\leq u_j \,\,\text{on}\,\, \Omega \,\,\text{and}\,\, v\leq g_j \,\,\text{on}\,\, \widetilde{\Omega} \}\big)^*.$$ 
Since $u_j\in C(\overline{\Om})$ and $g_j\in C(\overline{\O}),$ we infer that $\w_j\geq \min(\min\limits_{\Omega}u_j, \min\limits_{\tilde{\Omega}}g_j)$. Thus,  $\w_j\in SH_m^-(\O)\cap L^{\infty}(\O)\subset \E_m(\O).$
We put $$w_j=\begin{cases}
	\min(u_j,g_j)&\,\,\text{on}\,\,\Om\\
	g_j &\,\,\text{on}\,\,\O.
\end{cases}$$
Then we have $$\w_j=\big(\sup\{v\in SH^{-}_{m}(\tilde{\Omega}):  v\leq w_j \,\,\text{on}\,\, \widetilde{\Omega} \}\big)^*.$$
Lemma 4.2 in \cite{Pjmaa} implies that 
\begin{equation}\label{e3.2}
	H_m(\w_j)=0 \,\,\text{on}\,\, \{\w_j<w_j\}.
\end{equation}

\n On the other hand, it is easy to see that  and $\w_k\searrow M_{u,g}$ as $k\to\infty$
and
$$ \{M_{u,g}<s<\min(u,g)\}\cap\Om=\cup_{k=1}^{\infty}\big(\{\w_k<s<\min(u,g)\}\cap\Om\big),$$ so it is enough to prove that 
\begin{equation}\label{e3.1b} H_m(M_{u,g})=0\,\,\text{on}
\,\, \{\w_k<s<\min(u,g)\}\cap\Om.
\end{equation}
Indeed, from equality $\eqref{e3.2}$ we see that 
$$H_m(\w_j)=0 \,\,\text{on}\,\,\{ \w_j<s<w_j\}.$$ 

Let $k$ be an integer number. For $j\geq k$ we have
$$\big[\{u>s\}\cap\{g>s\}\cap \{\w_k<s\}\cap\Om\big]\subset \big[\{ \w_j<s<w_j\}\cap\Om\big].$$
Thus, we infer that
$$H_m(\w_j)=0 \,\,\text{on}\,\,\{u>s\}\cap\{g>s\}\cap \{\w_k<s\}\cap\Om.$$
It follows that
$$\max(u-s,0)\max(g-s,0)H_m(\w_j)=0 \,\,\text{on}\,\,\{\tilde{w}_k<s\}\cap\Om.$$

\n Letting $j\to\infty,$ since the equality 
\begin{align*}
	&2\max(u-s,0)\max(g-s,0)= \big( \max(u-s,0)+\max(g-s,0)\big)^2\\
	&-\big( \max(u-s,0)\big)^2 -\big( \max(g-s,0)\big)^2
\end{align*}
and Theorem  $\ref{th3.8HP17}$, we have
$$\max(u-s,0)\max(g-s,0)H_m(M_{u,g})=0 \,\,\text{on}\,\,\{\w_k<s\}\cap\Om.$$
According to Lemma 4.2 in $\cite{KH09}$, we deduce that
$$H_m(M_{u,g})=0\,\,\text{on}
\,\, \{\w_k<s<\min(u,g)\}\cap\Om.$$
That means we have equality $\eqref{e3.1b}$ and the proof of (ii) is complete.

\n (iii)
\n Firstly, we will prove that 
\begin{equation}\label{e3.3}
	H_m(\w_j)\leq H_m(u_j)  \,\,\text{on}\,\,\big(\{\w_j=u_j\}\cap \Om\big).
\end{equation}
Indeed, take a compact set $E\Subset \Omega\cap\{\w_j=u_j\}.$ Then for $\varepsilon >0$ we have $E\Subset \{\w_j+\varepsilon>u_j\}\cap\Omega.$ By Theorem 3.6 in $\cite{HP17}$ we have
\begin{align*}
	\int_{E}H_{m}(\w_j)&=\int_{E}\ind_{\{\w_j+\varepsilon>u_j\}}H_{m}(\w_j)\\
	&=\int_{E}\ind_{\{\w_j+\varepsilon>u_j\}}H_{m}(\max(\w_j+\varepsilon,u_j))\\
	&\leq \int_{E}H_{m}(\max(\w_j+\varepsilon,u_j)).
\end{align*}
Observe that $\Om$ is a open set and $\max(\w_j,u_j)=u_j$ on $\Om$, so we have $$H_{m}(\max(\w_j,u_j))=H_m(u_j)\,\,\text{on}\,\,\Om.$$  By Theorem $\ref{th3.8HP17}$, we get $H_{m}(\max(\w_j+\varepsilon,u_j))$ is weakly convergent to $H_{m}(u_j)$ when $\varepsilon\to 0.$ According to Lemma 1.9 in $\cite{De93}$ we get
$$\int_{E}H_{m}(\w_j) \le \varlimsup_{\ve \to 0} \int_{E}H_{m}(\max(\w_j+\varepsilon,u_j)) \le \int_E H_m(u_j).$$
This means that we have inequality $\eqref{e3.3}.$\\
\n Similarly, we also have 
\begin{equation}\label{e3.4}
	H_m(\w_j)\leq H_m(g_j)\,\,\text{on}\,\,\{\w_j=g_j\}.
\end{equation}
Combining equality $\eqref{e3.2}$, inequality $\eqref{e3.3}$ and inequality $\eqref{e3.4}$ we get
\begin{equation}\label{e3.5}
H_m(\w_j)\leq \ind_{\Om}H_m(u_j)+ H_m(g_j).
\end{equation}
By Theorem $\ref{th3.8HP17}$ when $j\to\infty$ we have $H_m(\w_j)\to H_m(M_{u,g})$ weakly on $\O,$ $H_m(g_j)\to H_m(g)$ weakly on $\O$ and $H_m(u_j)\to H_m(u)$ weakly on $\Om.$ From $\eqref{e3.5},$
let $j\to\infty$ we see that 
$$\ind_{\O\smallsetminus\partial\Om}H_m(M_{u,g})\leq \ind_{\O\smallsetminus\partial\Om}H_m(g)+\ind_{\Om}H_m(u).$$
Thus, it remains to prove that 
\begin{equation}\label{e3.6}
H_m(M_{u,g})\leq H_m(g)+\ind_{\Om}H_m(u) \,\,\text{on}\,\,\O\cap\partial\Om.
\end{equation}
Indeed, let a compact set $K\subset\O\cap\partial\Om.$ We need to prove that \begin{equation}\label{e3.6b}\int_{K}H_m(M_{u,g})\leq \int_{K}H_m(g).\end{equation} Fix $\ve>0$, it follows from the hypothesis $\int_{\Om}H_m(u)<+\infty$ that we can choose an open subset $G \Subset \Om$ such that $\int\limits_{\Om \setminus G} H_m (u)<\ve.$ Obviously, we have $K\subset\O\smallsetminus\overline{G}.$ Let $\chi\in C_{0}^{\infty}( \O\smallsetminus\overline{G}), 0\leq\chi\leq 1$ and $\chi=1 $ on $K.$ Note that $\ind_{K}$ is a upper semicontinuous function, so there exists a sequence $[\psi_j]\subset C(\O)\cap[0,1], \psi_j\searrow \ind_{K}$ on $\O$ as $j\to\infty.$
Now, by inequality $\eqref{e3.5}$ we have
\begin{align*}
\int_{K}H_m(M_{u,g})&\leq\int_{\O}\chi\psi_kH_m(M_{u,g})=\lim\limits_{j\to\infty}\int_{\O}\chi\psi_kH_m(\w_j)\\
&\leq \limsup\limits_{j\to\infty}\int_{\O}\chi\psi_kH_m(g_j)+\limsup\limits_{j\to\infty}\int_{\Om}\chi\psi_kH_m(u_j)\\
&\leq \int_{\O}\chi\psi_kH_m(g)+\limsup\limits_{j\to\infty}\int_{\Om\smallsetminus\overline{G}}\chi\psi_kH_m(u_j)\\
& (\text{by Lemma 1.9 in} \,\,\cite{De93}\,\, \text{with note that supp}H_m(g_j)\Subset\O.)\\
&\leq\int_{\O}\chi\psi_kH_m(g)+\limsup\limits_{j\to\infty}\int_{\Om\smallsetminus\overline{G}}H_m(u_j)\\
& \leq\int_{\O}\chi\psi_kH_m(g)+\limsup\limits_{j\to\infty}\int_{\Om}H_m(u_j)-\liminf\limits_{j\to\infty}\int_{\overline{G}}H_m(u_j)\\
&\leq\int_{\O}\chi\psi_kH_m(g)+\limsup\limits_{j\to\infty}\int_{\Om}H_m(u_j)-\liminf\limits_{j\to\infty}\int_{G}H_m(u_j)\\
&\leq\int_{\O}\chi\psi_kH_m(g)+\int_{\Om}H_m(u)-\liminf\limits_{j\to\infty}\int_{G}H_m(u_j)\\
& (\text{by Lemma 1.9 in} \,\,\cite{De93}\,\, \text{with note that supp}H_m(u_j)\Subset\Om.)\\
&\leq\int_{\O}\chi\psi_kH_m(g)+\int_{\Om}H_m(u)-\int_{G}H_m(u)\\
&\leq \int_{\O}\chi\psi_kH_m(g) + \int\limits_{\Om \setminus G} H_m (u)\\
& \leq \int_{\O}\psi_kH_m(g) +\ve.
\end{align*}
Letting $k\to\infty,$ by the Lebesgue's monotone convergence Theorem, we have 
$$\int_{K}H_m(M_{u,g})\leq \int_{K}H_m(g)+\ve. $$
Letting $\ve\to 0$ we obtain inequality $\eqref{e3.6b}.$ The proof of Part (A) of the Theorem is complete.

\n {\bf The proof of Part (B)}\\
(i) Repeating the same argument as in the proof of (i) of Part (A) with note that since $\varphi\in\mathcal{F}_m^{a}(\Omega)$ and $H_m(\hat{\varphi})\leq\ind_{\Omega}H_m(\varphi)$ we have $\hat{\varphi}\in\mathcal{F}_m^{a}(\tilde{\Omega}).$ Therefore, we get $M_{u,g}\in \mathcal{F}^a_{m}(g,\widetilde{\Omega}).$
\\ \n (ii) Using the same argument as in the proof of (ii) of Part (A) we have 
\begin{equation}\label{e1} H_m(M_{u,g})=0\,\,\text{ on}\,\, \{M_{u,g}<\min(u,g)\}\cap\Om.
\end{equation}
 Thus, it remains to prove that $$H_m(M_{u,g})=0\,\,\text{ on}\,\,\big((\O\smallsetminus K)\cap\{-\infty<M_{u,g}<g\}\big).$$
Indeed, the assumption $u\geq g$ on $\Om\smallsetminus K$ implies that $w_j\geq g$ 	on $\Om\smallsetminus K.$ Thus, we have $\big((\O\smallsetminus K)\cap\{\w_j<g\}\big)\subset\{\w_j<w_j\}.$
From equality $\eqref{e3.2}$ we infer that 
$$H_m(\w_j)=0 \,\,\text{ on}\,\,\big((\O\smallsetminus K)\cap\{\w_j<g\}\big).$$
Therefore, we have 	
$$H_m(\w_j)=0 \,\,\text{ on}\,\,\big((\O\smallsetminus K)\cap\{-\infty<\w_j<g\}\big).$$
Since
$$(\O\smallsetminus K)\cap\{-\infty<\w_j<g\}=\cup_{a,b\in\mathbb{Q}^-}(\O\smallsetminus K)\cap\{a<\w_j<b<g\} ,$$
we deduce that
$$H_m(\w_j)=0 \,\,\text{ on}\,\,(\O\smallsetminus K)\cap\{a<\w_j<b<g\}.$$
It follows that
$$\max(\w_j-a,0)\max(g-b,0)H_m(\w_j)=0\,\,\text{on}\,\, (\O\smallsetminus K)\cap\{\w_j<b\}.$$
Let $k$ be an integer number. For $j\geq k,$ since
$$(\O\smallsetminus K)\cap\{\w_k<b\}\subset (\O\smallsetminus K)\cap\{\w_j<b\},$$
we obtain
$$\max(\w_j-a,0)\max(g-b,0)H_m(\w_j)=0\,\,\text{on}\,\, (\O\smallsetminus K)\cap\{\w_k<b\}.$$
By Theorem 3.6 in $\cite{HP17}$ we have
$$\max(\w_j-a,0)\max(g-b,0)H_m(\max(\w_j,a))=0\,\,\text{on}\,\, (\O\smallsetminus K)\cap\{\w_k<b\}. $$
Letting $j\to\infty,$ by Lemma 3.8 in $\cite{HZ23}$ we get 
$$\max(M_{u,g}-a,0)\max(g-b,0)H_m(\max(M_{u,g},a))=0\,\,\text{on}\,\, (\O\smallsetminus K)\cap\{\w_k<b\}.$$
By Theorem 3.6 in $\cite{HP17}$ again, we infer that
$$\max(M_{u,g}-a,0)\max(g-b,0)H_m(M_{u,g})=0\,\,\text{on}\,\, (\O\smallsetminus K)\cap\{\w_k<b\}.$$
Since $\{M_{u,g}<b\}=\cup_{k=1}^{\infty}\{\w_k<b\},$ we have
$$\max(M_{u,g}-a,0)\max(g-b,0)H_m(M_{u,g})=0\,\,\text{on}\,\, (\O\smallsetminus K)\cap\{M_{u,g}<b\}.$$
It follows from Lemma 4.2 in $\cite{KH09}$ that
$$H_m(M_{u,g})=0\,\,\text{on}\,\, (\O\smallsetminus K)\cap\{a<M_{u,g}<b<g\}.$$ Thus, we get 
$$H_m(M_{u,g})=0\,\,\text{ on}\,\,\big((\O\smallsetminus K)\cap\{-\infty<M_{u,g}<g\}\big)$$
as desired.\\
\n (iii) We put 
$$w=\big(\sup\{v\in SH_m^-(\Om): v\leq u \,\,\text{on}\,\, K \}\big)^*.$$
It follows from $u\in\E_m(\Om)$ and local property of the class $\E_m(\Om)$ in $\F_m(\Om),$ we have $w\in\F_m(\Om).$  We will prove that $M_{u,g}=M_{w,g}$ on $\O.$ 
 To prove this, we have the following arguments:\\
a) Since $u\leq w,$ on $\Om$ we infer that $M_{u,g}\leq M_{w,g}.$\\
b) Assume that $SH_m^-(\O)\ni v\leq M_{w,g}.$\\
+) Since $v\leq w$ on $\Om$ we deduce that $v\leq u$ on $K.$ \\
++) On the other hand, from the hypothesis $u\geq g$ on $\Om\smallsetminus K$ and $v\leq g$ on $\O,$ we infer that $v\leq u$ on $\Om\smallsetminus K.$ \\
 Combining (+) and (++) we have $v\leq u$ on $\Om.$ By the definition of $M_{u,g}$ we have $v\leq M_{u,g}.$ This implies that $M_{w,g}\leq M_{u,g}.$\\
Combining (a) and (b) we get $M_{u,g}=M_{w,g}$ as desired.\\
\n Since $\O$ is an open set, so we have $H_m(M_{u,g})=H_m(M_{w,g})$ on $\O.$ It follows from $w\in\F_m(\Om)$ and Theorem 4.9 in \cite{T19} that $\int_{\Om}H_m(w)<+\infty.$ According to Part A of Theorem $\ref{t3.3}$ we get
\begin{equation*}H_m(M_{u,g})=H_m(M_{w,g})\leq \ind_{\Om}H_m(w)+ H_m(g) \,\,\text{on}\,\, \O.
\end{equation*}
\n  Note that $w$ is maximal on $\Om\smallsetminus K.$ It follows from Theorem 1.2 in \cite{Bl1} that  $H_m(w)=0$ on $\Om\smallsetminus K.$ 
This implies that 
\begin{equation}\label{eq3.11}H_m(M_{u,g})\leq \ind_{K}H_m(w)+ H_m(g) \,\,\text{on}\,\, \O.
\end{equation}
Repeating the argument as in the proof of inequality \ref{e3.3} we have 

\begin{equation}\label{eq3.7}H_m(M_{u,g})\leq H_m(u) \,\,\text{on}\,\, K\cap\{M_{u,g}=w\}\cap\{w>-\infty\} .
	\end{equation}
and
\begin{equation}\label{eq3.8}H_m(M_{u,g})\leq H_m(g) \,\,\text{on}\,\,\{M_{u,g}=g\}\cap\{g>-\infty\}.
\end{equation}
Since $u\leq w$ on $\Om,$ in view of Proposition $\ref{pro5.2HP17}$ we have
\begin{equation}\label{eq3.9}H_m(w)\leq H_m(u) \,\,\text{on}\,\, K\cap\{w=-\infty\}.
\end{equation}
\n Moreover, for all $m$-polar sets $A,$ it follows from $M_{u,g}\geq g+\hat{\varphi},$  Proposition $\ref{pro5.2HP17}$ and Lemma 5.6 in $\cite{HP17},$  we have
\begin{equation}\label{e3.12}
	\int_{A}H_m(M_{u,g})\leq \int_{A}H_m(g+\hat{\varphi})\leq \int_{A}H_m(g).
	\end{equation}
Combining the above argument we obtain that
\begin{equation}\label{e3.13}
H_m(M_{u,g})\leq \ind_{K\cap\{\u=u\}}H_m(u)+H_m(g).
\end{equation}
Indeed, it follows from inequality \ref{eq3.11} that inquality \ref{e3.13} is true on $\tilde{\Omega}\smallsetminus K.$ Moreover, since inequality \ref{eq3.9}, we only need to check inequality \ref{e3.13} on $K\cap\{w>-\infty\}.$ By inequality \ref{eq3.7} we see that inequality \ref{e3.13} works on $K\cap\{w>-\infty\}\cap \{ M_{u,g}=u\}.$ Thus, it remains to check inequality \ref{e3.13} on $K\cap\{w>-\infty\}\cap \{ M_{u,g}<u\}.$ According to inequality \ref{e1} we only need to check inequality \ref{e3.13} on $K\cap\{w>-\infty\}\cap \{ M_{u,g}<u\}\cap \{M_{u,g}=g\}.$ By inequality \ref{eq3.8} and inequality \ref{e3.12} we get desired.
The proof is complete.
\end{proof}

\section{Stability on $m$-capacity for maximal subextension of $m$-subharmonic functions with given boundary values}
First, we need the following Lemma.
\begin{lemma}\label{lm4.1}
Let $\Om\subset\O\Subset\mathbb{C}^n$ be $m$-hyperconvex domains and let $\{\Om_j\}$ be a sequence of bounded $m-$hyperconvex domains such that $\Om_j\Subset\Om_{j+1}\Subset\Om$ and $\Om=\cup_{j=1}^{\infty}\Om_j.$ Assume that $u\in\F_m^a(\Om)$ and consider $u_j$ defined by
$$u_j=\big(\sup\{v\in SH^{-}_m(\Om):v\leq u\,\,\text{on}\,\,C\Om_j\}\big)^*,$$
then $M_{u_j}\nearrow 0$ a.e. in $\O$ as $j\to\infty$.
\end{lemma}
\begin{proof}
Since $u_j\nearrow 0$ a.e. in $\Om,$ we also have $[M_{u_j}]$ is an increasing sequence. We put $w=(\sup\limits_{j\geq 1}M_{u_j})^*,$ then $M_{u_j}\nearrow w$ a.e. as $j\to+\infty.$ Since $u_j\geq u$ and $u\in\F_m(\Om),$ we obtain $u_j\in\F_m(\Om).$
In the view of Theorem $\ref{t3.3}$ we have $M_{u_j}\in\F_m(\O)$ and
\begin{equation}\label{e4.1}H_m(M_{u_j})\leq\ind_{\Om\cap\{M_{u_j}=u_j\}}H_m(u_j).
	\end{equation}
Since $\max(u_j,-1)\nearrow 0$ a.e. on $\Om$ and $H_m(u)$ puts no mass on all $m$-polar sets due to $u\in\F_m^a(\Om),$ according to Lemma $\ref{lm4.6PDjmaa}$ we have 
\begin{equation}\label{e4.2}\lim\limits_{j\to\infty}\int_{\Om}-\max(u_j,-1)H_m(u)=0.
	\end{equation}
Moreover, by Theorem $\ref{th3.8HP17},$ we have 
$$\max(w,-1)H_m(M_{u_j})\to\max(w,-1)H_m(w)$$
weakly as $j\to\infty$.\\
Now, by Lemma 1.9 in $\cite{De93}$ and inequality $\eqref{e4.1}$ we have 
\begin{align*}
	\int_{\O}-\max(w,-1)H_m(w)&\leq \liminf\int_{\O}-\max(w,-1)H_m(M_{u_j})\\
	&\leq \liminf\int_{\O}-\max(M_{u_j},-1)H_m(M_{u_j})\\
	&\leq\limsup\int_{\O}-\max(M_{u_j},-1)H_m(M_{u_j})\\
	&\leq \limsup\int_{\O}-\max(M_{u_j},-1)\ind_{\Om\cap\{M_{u_j}=u_j\}}H_m(u_j)\\
	&\leq \limsup\int_{\Om}-\max(u_j,-1)H_m(u_j)
\end{align*}
Since $u_j\geq u$ on $\Om,$ it follows from Lemma $\ref{lm2.7Gasmi}$ and  equality $\eqref{e4.2}$ that
$$\int_{\Om}-\max(u_j,-1)H_m(u_j)\leq \int_{\Om}-\max(u_j,-1)H_m(u) =0 .$$
Therefore, we have $\int\limits_{\O}-\max(w,-1)H_m(w)\leq 0.$ This implies that $$\mu= -\max(w,-1)H_m(w)=0\,\,\text{ on}\,\, \O.$$
Note that, since $M_{u_j}\in\F_m(\O),$ we also have $w\in\F_m(\O).$ We claim that 
$$H_m(w)=0 \,\,\text{ on}\,\, \O.$$ 
Indeed, let $\chi$ be a continuous function with compact support on $\O.$ Then we have
\begin{align*}
\int_{\O}\chi H_m(w)=&\int_{\O\cap\text{supp}\chi}\chi H_m(w)\\
&=\int_{\O\cap\text{supp}\chi}\dfrac{\chi}{-\max(w,-1)} [-\max(w,-1)]H_m(w)\\
&=\int_{\O\cap\text{supp}\chi}\dfrac{\chi}{-\max(w,-1)}d\mu\\
&=0.
\end{align*}
This proves the claim. According to Remark 5.4 in $\cite{HP17},$ we get $w=0$ on $\O.$ Thus, the proof of Lemma is complete.
\end{proof}

\n Next, we state a version of "comparison principle" for the class $\mathcal{N}_m(f,\Om).$
\begin{lemma}\label{t4.2}
Let $f\in\E_m(\Om).$ Assume that $u\in\mathcal{N}_m^a(f,\Om)$ and $v\in\E_m(\Om), v\leq f$ are such that\\
\n (i) $H_m(u)=0$ on $\{-\infty<u<v\}.$\\
\n (ii) $\int_{\Om}H_m(u)<+\infty.$\\
Then we have $u\geq v$ on $\Om.$
\end{lemma}
\begin{proof}
We put $w=\max(u,v).$ It follows that $w\in\mathcal{N}_m^a(f,\Om)$ and $u\leq w$ on $\Om.$
We have $ \{-\infty<u<v\}=\{-\infty<u<w\}.$ By (i) we have
\begin{equation}\label{eq4.3}
	H_m(u)=0 \,\,\text{on}\,\, \{-\infty<u<w\}.
\end{equation}
Repeating the argument as in inequality \ref{e3.3} we obtain
\begin{equation}\label{eq4.4}
H_m(u)\leq H_m(w) \,\,\text{on}\,\,\{u=w\}\cap\{u>-\infty\}
\end{equation}
On the other hand, by hypothesis (ii) and Proposition 3.5 in $\cite{HZ23},$ we have
\begin{equation}\label{eq4.5} \ind_{\{f=-\infty\}}H_m(f)=\ind_{\{u=-\infty\}}H_m(u).
	\end{equation}
Moreover, since $w\geq u$ and $\int_{\Om}H_m(u)<+\infty,$ according to Lemma $\ref{lm2.7Gasmi},$ we have
$$\int_{\Om}H_m(w)\leq \int_{\Om}H_m(u)<+\infty.$$
Therefore, by Proposition 3.5 in $\cite{HZ23}$ again, we have
\begin{equation}\label{eq4.6}\ind_{\{f=-\infty\}}H_m(f)=\ind_{\{w=-\infty\}}H_m(w).
	\end{equation}
Combining equality $\eqref{eq4.5}$ and equality $\eqref{eq4.6}$, we obtain
\begin{equation}\label{eq4.7}
\ind_{\{u=-\infty\}}H_m(u)=\ind_{\{w=-\infty\}}H_m(w)\leq \ind_{\{u=-\infty\}}H_m(w) .
\end{equation}
Combining equality $\eqref{eq4.3}$, inequality $\eqref{eq4.4} $ and inequality $\eqref{eq4.7},$ we infer that
$$H_m(u)\leq H_m(w)\,\,\text{on}\,\,\Om.$$
By Proposition 3.6 in $\cite{HZ23},$ we imply that $u\geq w$ on $\Om.$ Thus, we infer that $u=w$ on $\Om.$ Therefore, we have $u\geq v$ on $\Om.$ The proof of Lemma is complete.
\end{proof}

\n Now, we prove main result in this Section.
\begin{theorem}
	Let $\Om\subset\O\Subset\mathbb{C}^n$ be $m-$hyperconvex domains and $f\in\E_m(\Om), g\in\E_m(\O), w\in\F_m^a(f,\Om)$ are such that $f\geq g$ on $\Om$ and
	\begin{equation}\label{eq4.8}\int_{\Om}H_m(w)+\int_{\O}H_m(g)<+\infty.
	\end{equation}
	Assume that $[u_j]\subset\F_m^{a}(f,\Om), u_j\geq w$ and $u_j$ converges to $u\in SH_m(\Omega)$ in $m-$ capacity on $\Om.$ Then $M_{u_j,g}$ converges to $M_{u,g}$ in $m$-capacity on $\O.$ 
\end{theorem}
\begin{proof}
	It is enough to prove that there exist $M_{u,g}\in SH_m(\O)$ such that from any subsequence of the given sequence $\{M_{u_{j},g}\}$ one can extract sub-subsequence $\{M_{u_{j_k},g}\}$ converging to $M_{u,g}$ in $m$-capacity on $\O$.
	Obviously, it follows from the constructions of $M_{w,g}$ and $M_{u_{j},g}$ that we have $M_{w,g}\leq M_{u_{j},g}.$ Thus, by Lemma 2.5 in \cite{AAG20}, choose $[M_{u_{j_k},g}]$ from any subsequence such that it converges weakly to an $m$-subharmonic function $\phi\in SH_m(\O)$ in $L^1_{\loc}(\O)$. It remains to prove that $M_{u_{j_k},g}$ converges to $\phi$ in $m$-capacity on $\O$ and $\phi=M_{u,g}$ on $\O$.
	Indeed, replacing $[M_{u_{j_k},g}]$ by a subsequence, we can assume that $M_{u_{j_k},g}\to\phi$ a.e. on $\O.$  From now to the end of this proof, to make it convenient for presentation without changing argument, we will use symbol $M_{u_j,g}$ instead of $M_{u_{j_k},g}$. It means that we will prove that $M_{u_j,g}$ converges to $\phi$ in $m$-capacity on domain $\O$ with hypothesis $M_{u_j,g}\to \phi$ a.e. on domain $\O$. \\
	Taking $\phi_k=\big(\sup\limits_{j\geq k} M_{u_{j},g}\big)^* \,\,\text{on}\,\,\O$, to prove the convergence of $M_{u_j,g}\to\phi$ in $m$-capacity, Theorem A in \cite{HZ23} says that it suffices to prove the following equation
	
	\begin{equation}\label{e4.9}\lim\limits_{k\to\infty}\int_{\O}\big[\max( \dfrac{\phi_k}{r},\rho)-\max(\dfrac{M_{u_{k},g}}{r},\rho)\big]H_m(M_{u_{k},g})=0 \,\text{for all}\, r>0.
	\end{equation}
	We split the proof into two cases:\\
	{\bf Case 1.} Assume that $w\geq g$ on $\Om\smallsetminus K$ for some compact subset $K$ of $\Om.$ Since $u_j\geq w\geq g$ on $\Om\smallsetminus K,$ according to Theorem $\ref{t3.3}$ we have
	\begin{equation}\label{eq4.9}
		H_m(M_{u_j,g})\leq \ind_{K\cap\{M_{u_j,g}=u_j\}}H_m(u_j)+H_m(g)\,\,\text{on}\,\,\O.
	\end{equation}
	
	\n We put $\psi_k=\big(\sup\limits_{j\geq k} u_{j}\big)^* \,\,\text{on}\,\,\Om$ and take $h\in\E_m^0(\Om)\cap C(\Om)$ and $\rho\in\E_m^0(\O)\cap C(\O)$ such that $K\subset \Om\cap \{h=\rho\}.$\\
	Since $M_{u_{k},g}\leq \phi_k\leq g$ on $\O$ so for some real number $r>0$ we have
	$$\max(\dfrac{\phi_k}{r},\rho)=\max(\dfrac{M_{u_{k},g}}{r},\rho)=\rho \,\,\text{on}\,\,\O\cap\{g=-\infty\}.$$
	Combining the above equality and inequality $\eqref{eq4.9}$ we have
	\begin{equation}
		\begin{split}\label{eq4.10}
			&\limsup\limits_{k\to\infty}\int_{\O}\big[\max( \dfrac{\phi_k}{r},\rho)-\max(\dfrac{M_{u_{k},g}}{r},\rho)\big]H_m(M_{u_{k},g})\\
			&\leq \limsup\limits_{k\to\infty}\int_{K\cap\{M_{u_{k},g}=u_{k} \}}\big[\max( \dfrac{\phi_k}{r},\rho)-\max(\dfrac{M_{u_{k},g}}{r},\rho)\big]H_m(u_{k})\\
			&+\limsup\limits_{k\to\infty}\int_{\O}\big[\max( \dfrac{\phi_k}{r},\rho)-\max(\dfrac{M_{u_{k},g}}{r},\rho)\big]H_m(g)\\
			&\leq \limsup\limits_{k\to\infty}\int_{\Om}\big[\max( \dfrac{\psi_k}{r},h)-\max(\dfrac{u_{k}}{r},h)\big]H_m(u_{k})\\
			&+\limsup\limits_{k\to\infty}\int_{\O\cap\{g>-\infty\}}\big[\max( \dfrac{\phi_k}{r},\rho)-\max(\dfrac{M_{u_{k},g}}{r},\rho)\big]H_m(g).
		\end{split}
	\end{equation}
	It follows from $M_{u_{k},g}\to\phi$ a.e. on $\O$ that  $\max(\dfrac{M_{u_{k},g}}{r},\rho)\to \max(\dfrac{\phi}{r},\rho)$ a.e. on $\O$. Moreover, note that $\ind_{\{g>-\infty\}}H_m(g)$ is a measure on $\O$ that vanishes on all $m-$ polar sets and \begin{align*}\int_{\O}-\max(\dfrac{M_{u_{k},g}}{r},\rho)\ind_{\{g>-\infty\}}H_m(g)&\leq \int_{\O}-\rho H_m(g)\\
		&\leq -\min\limits_{\O}\rho\int_{\O} H_m(g)<+\infty.
	\end{align*}
	Thus, by Lemma $\ref{lm4.6PDjmaa}$ we have 
	\begin{equation}\label{eq4.11}
		\lim\limits_{k\to\infty}\int_{\O}\max(\dfrac{M_{u_{k},g}}{r},\rho)\ind_{\{g>-\infty\}}H_m(g)=\int_{\O}\max(\dfrac{\phi}{r},\rho)\ind_{\{g>-\infty\}}H_m(g).
	\end{equation}
	Obvious, we have $\phi_k\searrow\phi$ on $\O.$ Repeating the above argument we have
	\begin{equation}\label{eq4.12}
		\lim\limits_{k\to\infty}\int_{\O}\max(\dfrac{\phi_k}{r},\rho)\ind_{\{g>-\infty\}}H_m(g)=\int_{\O}\max(\dfrac{\phi}{r},\rho)\ind_{\{g>-\infty\}}H_m(g).
	\end{equation}
	Combining equality $\eqref{eq4.11}$ and equality $\eqref{eq4.12}$ we get
	\begin{equation}\label{eq4.13}
		\limsup\limits_{k\to\infty}\int_{\O\cap\{g>-\infty\}}\big[\max( \dfrac{\phi_k}{r},\rho)-\max(\dfrac{M_{u_{k},g}}{r},\rho)\big]H_m(g)=0.
	\end{equation}
	Moreover, it follows from $u_j$ converging to $u$ in $m$-capacity on $\Om$ and Theorem A in $\cite{HZ23}$  that
	\begin{equation}\label{eq4.14}
		\limsup\limits_{k\to\infty}\int_{\Om}\big[\max( \dfrac{\psi_k}{r},h)-\max(\dfrac{u_{k}}{r},h)\big]H_m(u_{k})=0
	\end{equation}
	Combining inequality $\eqref{eq4.10}$, equality $\eqref{eq4.13}$ and equality $\eqref{eq4.14},$ we infer that 
	$$\lim\limits_{k\to\infty}\int_{\O}\big[\max( \dfrac{\phi_k}{r},\rho)-\max(\dfrac{M_{u_{k},g}}{r},\rho)\big]H_m(M_{u_{k},g})=0. $$
	Thus, we obtain equality \ref{e4.9}. Therefore, we infer that $M_{u_{k},g}$ converges to $\phi$ in $m$-capacity on $\O$ as desired.
	
	Now we will prove that $M_{u,g}=\phi$ on $\O.$ Obviously, we have $\phi_k \leq \psi_k$ on $\Omega$. Letting $k\to\infty$, we get $\phi\leq u$ on $\Omega$. This implies that $\phi\leq M_{u,g}$. We will use Lemma $\ref{t4.2}$ to prove that $\phi\geq M_{u,g}$. Indeed, by $u_{j}\geq w$ we infer that $M_{u_j,g}\geq M_{w,g}\in\F_m^a(g,\O).$ Since  $M_{u_j,g}\to\phi$ a.e. on $\O,$ we deduce that $\phi\geq M_{w,g}$ a.e. on $\O.$ By the subharmonicity, we get $\phi\geq M_{w,g}$ on $\O.$ Thus, we have $\phi\in \F^a_m(g,\O).$\\
	\n Moreover, since $u_j$ converges to $u$ in $m-$ capacity on $\Om,$ we infer that $u_j$ converges weakly to $u$ on $\Om.$ This implies that there exists a subsequence $u_{j_k}$ of $[u_j]$ such that $u_{j_k}$ converges to $u$ a.e. on $\Om.$ Since $u_{j_k}\geq w$ on $\Om,$ we implies that $u\geq w$ a.e. on $\Om.$ By the subharmonicity, we get $u\geq w$ on $\Om.$ Since $w\in \F_m^a(f,\Om)$, we deduce that $u\in \F_m^a(f,\Om).$ By Theorem $\ref{t3.3},$ we get $M_{u,g}\in\F_m^a(g,\O).$ \\
	\n We claim that $\int_{\O}H_m(\phi)<+\infty$. Indeed, according to Theorem $\ref{th3.8HP17},$ we have $H_m(M_{u_{j},g})$ converges weakly to $H_m(\phi)$ as $j\to\infty.$ By Lemma 1.9 in $\cite{De93}$, inequality $\eqref{eq4.9}$ and the hypothesis $\eqref{eq4.8}$,  we have
	\begin{align*}
		\int_{\O}H_m(\phi)&\leq \liminf_{j\to\infty}\int_{\O}H_m(M_{u_{j},g})\\
		&\leq \liminf_{j\to\infty}\int_{\O}\big[\ind_{\Om}H_m(u_{j})+H_m(g)\big]\\
		& \leq \liminf_{j\to\infty}\int_{\Om}H_m(u_{j})+\int_{\O}H_m(g)\\
		&\leq \int_{\Om}H_m(w) + \int_{\O}H_m(g)\\
		&< \infty,
	\end{align*}
	where the fourth inequality is due to $u_{j}\geq w$ and  Lemma $\ref{lm2.7Gasmi}.$	\\
	\n Thus, to apply Lemma $\ref{t4.2}$ to get $\phi\geq M_{u,g}$, it remains to prove that
\begin{equation}\label{e4.16}	H_m(\phi)=0\,\,\text{on}\,\,\{-\infty<\phi<M_{u,g}\}.
	\end{equation}
	
	\n Note that, by Theorem $\ref{t3.3}$ we have $H_m(M_{u_{j},g})=0 $ on
	\begin{equation}\label{eq4.15}
		\big((\O\smallsetminus K)\cap\{-\infty<M_{u_{j},g}<g\}\big)\cup\big(\Om\cap\{-\infty< M_{u_{j},g}<\min(u_{j},g)\}\big).
	\end{equation}
	\n Using the assumption $M_{u_j,g}$ converges to $\phi$ in $m$-capacity, we will prove that
	\begin{equation}\label{eq4.16}
		H_m(\phi)=0 \,\,\text{on}\,\,\big((\O\smallsetminus K)\cap\{-\infty<\phi<g\}\big)\cup\big(\Om\cap\{-\infty<\phi<\min(u,g)\}\big).
	\end{equation}
	Indeed, we have
	$$(\O\smallsetminus K)\cap\{-\infty<M_{u_j,g}<g\}=\cup_{a,b\in\mathbb{Q}^-}\big((\O\smallsetminus K)\cap\{a<M_{u_j,g}<b<g\}\big)$$
	and $$ \Om\cap\{-\infty<M_{u_j,g}<\min(u_{j},g)\}=\cup_{c,d\in\mathbb{Q}^-}\big(\Om\cap\{c<M_{u_j,g}<d<\min(u_{j},g)\}\big).$$
	It follows from the equality $\eqref{eq4.15}$ that $H_m(M_{u_j,g})=0 $ on 
	$$ \big((\O\smallsetminus K)\cap\{a<M_{u_j,g}<b<g\}\big)\cup \big(\Om\cap\{c<M_{u_j,g}<d<\min(u_{j},g)\}\big).$$
	This implies that 
	$$\max(M_{u_j,g}-a,0)\max(g-b,0)H_m(M_{u_j,g})=0\,\,\text{on}\,\,\big((\O\smallsetminus K)\cap\{M_{u_j,g}<b\}\big) $$
	and 
	$$ \max(M_{u_j,g}-c,0)\max(u_j-d,0)\max(g-d,0)H_m(M_{u_j,g})=0\,\,\text{on}\,\,\big(\Om\cap\{M_{u_j,g}<d\}\big).$$
	For an integer number $j\geq k$ we have
	$$\big((\O\smallsetminus K)\cap\{\phi_{k}<b\}\big)\subset \big((\O\smallsetminus K)\cap\{\phi_{j}<b\}\big)\subset  \big(\O\smallsetminus K)\cap\{M_{u_j,g}<b\}\big) $$
	and $$\big( \Om\cap\{\phi_{k}<d\}\big)\subset\big( \Om\cap\{\phi_{j}<d\}\big)\subset \big(\Om\cap\{M_{u_j,g}<d\}\big).$$
	Therefore, we have 
	$$\max(M_{u_j,g}-a,0)\max(g-b,0)H_m(M_{u_j,g})=0\,\,\text{on}\,\,\big((\O\smallsetminus K)\cap\{\phi_{k}<b\}\big) $$
	and 
	$$ \max(M_{u_j,g}-c,0)\max(u_j-d,0)\max(g-d,0)H_m(M_{u_j,g})=0\,\,\text{on}\,\,\big(\Om\cap\{\phi_{k}<d\}\big).$$
	\n By Theorem 3.6 in $\cite{HP17}$ we have
	$$\max(M_{u_j,g}-a,0)\max(g-b,0)H_m(\max(M_{u_j,g};a))=0\,\,\text{on}\,\,\big((\O\smallsetminus K)\cap\{\phi_{k}<b\}\big) $$
	and 
	$$ \max(M_{u_j,g}-c,0)\max(u_j-d,0)\max(g-d,0)H_m(\max(M_{u_j,g};c))=0$$ on $$\big(\Om\cap\{\phi_{k}<d\}\big).$$
	Letting $j\to\infty,$ by Lemma 3.8 in $\cite{HZ23}$ we get 
	$$\max(\phi-a,0)\max(g-b,0)H_m(\max(\phi,a))=0\,\,\text{on}\,\,\big((\O\smallsetminus K)\cap\{\phi_{k}<b\}\big) $$
	and 
	$$ \max(\phi-c,0)\max(u-d,0)\max(g-d,0)H_m(\max(\phi,d))=0\,\,\text{on}\,\,\big(\Om\cap\{\phi_{k}<d\}\big).$$
	By Theorem 3.6 in $\cite{HP17}$ again, we infer that
	$$\max(\phi-a,0)\max(g-b,0)H_m(\phi)=0\,\,\text{on}\,\,\big((\O\smallsetminus K)\cap\{\phi_{k}<b\}\big) $$
	and 
	$$ \max(\phi-c,0)\max(u-d,0)\max(g-d,0)H_m(\phi)=0\,\,\text{on}\,\,\big(\Om\cap\{\phi_{k}<d\}\big).$$
	Since $\{\phi<b\}=\cup_{k=1}^{\infty}\{\phi_k<b\}$ and $\{\phi<d\}=\cup_{k=1}^{\infty}\{\phi_k<d\},$ we have
	$$\max(\phi-a,0)\max(g-b,0)H_m(\phi)=0\,\,\text{on}\,\,\big((\O\smallsetminus K)\cap\{\phi<b\}\big) $$
	and 
	$$ \max(\phi-c,0)\max(u-d,0)\max(g-d,0)H_m(\phi)=0\,\,\text{on}\,\,\big(\Om\cap\{\phi<d\}\big).$$
	It follows from Lemma 4.2 in $\cite{KH09}$ that
	$$H_m(\phi)=0\,\,\text{on}\,\, (\O\smallsetminus K)\cap\{a<\phi<b<g\}$$
	and
	$$H_m(\phi)=0\,\,\text{on}\,\,\Om\cap\{c<\phi<d<\min(u,g)\}.$$
	Thus, we get equality $\eqref{eq4.16}$ as  desired.\\
	Obviously, we have
	$$\{-\infty<\phi<M_{u,g}\}\subset \big((\O\smallsetminus K)\cap\{-\infty<\phi<g\}\big)\cup\big(\Om\cap\{-\infty<\phi<\min(u,g)\}\big).$$
	Thus, it follows from equality $\eqref{eq4.16}$ that
	$$H_m(\phi)=0\,\,\text{on}\,\,\{-\infty<\phi<M_{u,g}\}.$$
	It means that we get equality \ref{e4.16} and we finish the proof of Case 1.\\
	
	\n {\bf Case 2.} General Case.\\
	\n Let $\{\Om_k\}$ be a sequence of bounded $m-$ hyperconvex domains such that $\Om_k\subset\Om_{k+1}$ and $\Om=\cup_{k=1}^{\infty}\Om_k.$ Since $w\in\F_m^a(f,\Om),$ there exists a function $\xi\in\F_m^a(\Om) $ such that $f\geq w\geq f+\xi.$
	We consider
	$$\xi_k=\big(\sup\{v\in SH_m^-(\Om): v\leq\xi\,\,\text{on}\,\,\Om\smallsetminus\Om_k\}\big)^*,$$
	then we get $\xi_k\nearrow 0$ a.e. on $\Om.$ We set $\eta_k=M_{\xi_k}$ then by Lemma \ref{lm4.1} we have $\eta_k\nearrow 0$ a.e on $\O$.\\
	 Now, on $\Om\smallsetminus \Om_k$ we have
	$$u_j\geq w\geq f+\xi\geq f+\xi_k\geq g+\eta_k.$$
	By Case 1, we have 
	$M_{u_{j},g+\eta_k}\to M_{u,g+\eta_k}$
	in $m-$capacity on $\O$ as $j\to\infty.$\\
	\n Obviously, we have 
	$ M_{u_{j},g}\geq M_{u_{j},g+\eta_k}$ and
	$ M_{u,g}\geq M_{u,g+\eta_k}.$\\
	\n On the other hand, it is easy to see that 
	$$ M_{u_{j},g}+\eta_k\leq M_{u_{j},g+\eta_k}$$
	and 
	$$ M_{u,g}+\eta_k\leq M_{u,g+\eta_k}.$$
	Thus, we have
	\begin{align*}
		|M_{u_{j},g}-M_{u,g}|&\leq |M_{u_{j},g}-M_{u_{j},g+\eta_k}|+|M_{u_{j},g+\eta_k}-M_{u,g+\eta_k}|\\
		&					    +|M_{u,g+\eta_k}-M_{u,g}|\\
		&\leq -\eta_k+ |M_{u_{j},g+\eta_k}-M_{u,g+\eta_k}|-\eta_k\\
		&\leq |M_{u_{j},g+\eta_k}-M_{u,g+\eta_k}| -2 \eta_k
	\end{align*}
	This implies that $M_{u_{j_k},g}\to M_{u,g}$ in $m-$ capacity on $\O.$
	The proof of the Theorem is complete.
\end{proof}
\section*{Declarations}
\subsection*{Ethical Approval}
This declaration is not applicable.
\subsection*{Competing interests}
The authors have no conflicts of interest to declare that are relevant to the content of this article.
\subsection*{Funding }
No funding was received for conducting this study.
\subsection*{Availability of data and materials}
This declaration is not applicable.

\end{document}